\documentclass{amsart}

\usepackage{graphicx}
\usepackage{amsmath, amssymb, amsthm}
\usepackage[mathcal]{eucal}
\usepackage{enumerate}
\usepackage[dvipsinames]{xcolor}
\definecolor{midnightblue}{HTML}{0059b3}
\definecolor{chromered}{HTML}{f14233}
 \RequirePackage[colorlinks,citecolor=midnightblue,urlcolor=midnightblue,breaklinks=true,
linkbordercolor=midnightblue,linkcolor=midnightblue]{hyperref}

\newtheorem{theorem}{Theorem}
\newtheorem{conjecture}[theorem]{Conjecture}
\newtheorem{lemma}{Lemma}
\newtheorem{claim}{Claim}

\newtheorem{corollary}[theorem]{Corollary}
\newtheorem{OldTheorem}{Theorem}

\def\Id{{\rm Id}}
\def\Proj{{\rm Proj}}

\def\dist{{\rm dist}}

\def\sgn{{\rm sgn}}

\def\ti{\tilde}

\def\ZZ{\ensuremath{\mathbb Z}}

\def\D{\ensuremath{\mathbb D}}
\def\T{\ensuremath{\mathbb T}}
\def\R{\ensuremath{\mathbb R}}
\def\Z1{\ensuremath{\mathbf 1}}

\numberwithin{equation}{section}

\title[Wavelet Malmquist-Takenaka system]{Almost everywhere convergence of a wavelet-type Malmquist-Takenaka series}
\author{Gevorg Mnatsakanyan}
\address{Institute of Mathematics of the National Academy of Sciences of Armenia, Marshal Baghramyan Ave. 24/5, Yerevan 0019, Armenia}
\curraddr{}
\email{mnatsakanyan\_g@yahoo.com}
\date{December 2025}

\begin{document}

\begin{abstract}
    The Malmquist-Takenaka (MT) system is a complete orthonormal system in $H^2(\T)$ generated by an arbitrary sequence of points $a_n$ in the unit disk with $\sum_n (1-|a_n|) = \infty$. The point $a_n$ is responsible for multiplying the $n$th and subsequent terms of the system by a M\"obius transform taking $a_n$ to $0$. One can recover the classical trigonometric system, its perturbations or conformal transformations, as particular examples of the MT system.
    However, for many interesting choices of the sequence $a_n$, the MT system is less understood. In this paper, we consider a wavelet-type MT system and prove its almost everywhere convergence in $H^2(\T)$.
\end{abstract}
\maketitle

\section{Introduction}

Let $(a_n)_{n=1}^{\infty}$ be a sequence of points inside the unit disk such that
\begin{equation}\label{completeness condition}
\sum_{n=1}^{\infty} ( 1 - |a_n| ) = +\infty.
\end{equation}

The associated Blaschke products and Malmquist-Takenaka (MT) basis \cite{Malmquist, Takenaka} are defined as
\begin{equation}\label{BlaschkeMT}
B_n (z) = \prod_{j=1}^n \frac{\bar a_j}{|a_j|} \frac{z - a_j}{1 - \overline{a_j} z }, \quad\quad \phi_n (z) = B_{n} (z) \frac{ \sqrt{1 - | a_{n+1} |^2 } }{ 1 - \overline{ a_{n+1} } z } \, ,
\end{equation}
for $n\geq 0$, where we put $B_0=1$. It was shown in \cite{Coif2} that $(\phi_n)_{n=0}^\infty$ is an orthonormal system in $H^p(\D)$, $1<p<\infty$. The condition \eqref{completeness condition} ensures the completeness of the system. If $a_n\equiv 0$, the corresponding MT system reduces to the classical trigonometric system, that is $\phi_n(z)=z^n$. Other choices of the sequence $a_n$ enable various scenarios for the behavior of the system.

The Malmquist-Takenaka system has been studied in the context of information theory \cite{Qian,Qian2} where it is sometimes called the adaptive Fourier transform.
As the name suggests, one could try to fix a function $f\in H^2(\D)$ and choose a sequence $a_n$ depending on $f$ in a way that would improve the rate of convergence of the series
\begin{equation}\label{mt partial sum}
    \sum_{j=0}^n \langle f, \phi_j \rangle \phi_j(e^{i\theta})
\end{equation}
compared, say, to the classical Fourier series. One strategy towards that goal could be to choose $a_n$ by a greedy algorithm maximizing the magnitude of the coefficient $\langle f, \phi_n \rangle$ at the $n$th step. This has been explored in \cite{Coif2}. Another approach would be to make as many of the coefficients equal to $0$ as possible and to look at an appropriate subsequence of the partial sums. This is accomplished by the nonlinear phase unwinding decomposition, which can equivalently and perhaps more naturally be defined through the Blaschke factorization. Taking $F\in H^2(\D)$ we define $B_1$ to be the Blaschke product having the same zeros as $F-F(0)$ and put
$$F_1 := \frac{F-F(0)}{B_1} \, .$$
Iterating, let $B_n$ be the Blaschke product having the same zeros as $F_{n-1}-F_{n-1}(0)$ and
$$
F_n := \frac{F_{n-1}-F_{n-1}(0)}{B_n} \, .
$$
We arrive at the formal series called, the nonlinear phase unwinding,
\begin{equation}\label{unwinding}
    F(z) = F(0)+F_1(0)B_1(z)+\dots + F_n(0)B_1(z)\dots B_n(z)+\cdots \, .
\end{equation}
The unwinding series \eqref{unwinding} coincides with the MT series \eqref{mt partial sum} if we choose $a_n$ to be the zeros of the Blaschke products $B_1, B_2$ and so on. The unwinding series was introduced in \cite{Nah} with numerical simulations suggesting a very fast convergence rate. A variety of interesting convergence results are obtained in \cite{Coif1}. For more results and discussion we refer to \cite{Steinerberger1,steinerberger2, Coif3}. Whether the unwinding series \eqref{unwinding} converges almost everywhere for $F\in H^2(\D)$ or whether \eqref{unwinding} converges "fast" for "most" functions remain interesting open questions.

In the present paper we investigate the almost everywhere convergence property of the MT series for two particular sequences.
By standard techniques a positive result on the almost everywhere convergence follows from the boundedness of maximal partial sum operator.
Let us introduce the maximal partial sum operator of the MT series. For an arbitrary sequence $\mathbf{a} = (a_n)_{n=1}^\infty$ and $f\in L^2(\T)$, we denote
\begin{equation}\label{carleson op intro}
    T^{\mathbf{a}} f (x) := \sup_{n \geq 0} \left|  \sum_{j=0}^n \langle f, \phi_j \rangle \phi_j(e^{ix}) \right| \, .
\end{equation}

Throughout the paper $[x]$ denotes the integer part of $x$. Our main theorem is the following.
\begin{theorem}\label{main}
    For any $\frac{1}{2}<r<1$ let $a_n = r e^{2\pi in(1-r)}$, $1\leq n \leq \left[\frac{1}{1-r} \right]$. Consider the sequence $\mathbf{a}_r = (a_n)_{n=1}^{[1/(1-r)]}$. Then, for any $f\in L^2(\mathbb{T})$
    \begin{equation}\label{maininequality}
        \| T^{\mathbf{a_r}} f \|_{L^2(\T)} \lesssim \| f\|_{L^2(\T)}\, ,
    \end{equation}
    with the implicit constant independent of $r$.
\end{theorem}
We have the following corollary.
\begin{corollary}\label{aeconvergence}
    For $n\in \mathbb{N}$, let $b_n := (1-2^{-[\log_2 n]}) e^{2\pi i n 2^{-[\log_2 n ]}}$ and $\mathbf{b} = (b_n)_{n=1}^{\infty}$. Then, for any $f\in L^2(\T)$,
    \begin{equation}
        \| T^{\mathbf{b}} f \|_{L^2(\T)} \lesssim \|f\|_{L^2(\T)} \, . 
    \end{equation}
    In particular, MT series \eqref{mt partial sum} associated to $\mathbf{b}$ converges almost everywhere on $\T$.
\end{corollary}
Letting $r_m = 1-2^m$, $m\geq 1$, the sequence $\mathbf{b}$ can be obtained by concatenating $\mathbf{a}_{r_1}, \mathbf{a}_{r_2}$ and so on. 

We call the MT systems resulting from $\mathbf{a_r}$ and $\mathbf{b}$ wavelet-type MT systems. Similar systems were introduced and studied in \cite{pap,Feichtinger2013, papschipp} in the context of multisresolution of Hardy spaces and are called hyperbolic wavelet systems. We remark that our proofs also work for sequences of the form $a_n =(1-r) e^{i Dn (1-r)}$ with some parameter $D$ and the implicit constant in Theorem \ref{main} depending on $D$. In particular, this includes the MT systems considered in \cite[Formulas (2.6) and (2.7)]{pap}. However, for the sake of keeping the exposition simple and emphasizing the method of the proof we prefer to restrict our attention to $\mathbf{a_r}$ and $\mathbf{b}$.

Other almost everywhere convergence results for the MT series were previously obtained by the author in \cite{MNATSAKANYAN2022109461}. We want to compare Theorem \ref{main} to
\begin{OldTheorem}[\cite{MNATSAKANYAN2022109461}]\label{oldmain}
    Let $0<r<1$ and let $\mathbf{c} = (c_n)_{n=1}^\infty$ be an arbitrary sequence such that $|c_n|\leq r$ for all $n$, then
    \begin{equation}\label{upperbound compact}
        \|T^{\mathbf{c}} \|_{L^2(\T)\to L^2(\T)} \lesssim \sqrt{\log \frac{1}{1-r}} \, .
    \end{equation}
   Moreover, for $d_n = re^{2\pi i n (1-r)\log \frac{1}{1-r}}$, $1\leq n \leq \frac{1}{(1-r)\log \frac{1}{1-r}}$, and $\mathbf{d_r} = (d_n)_{n=1}^{[1/(1-r)\log \frac{1}{1-r}]}$, there exists $f\in L^2(\T)$ such that
    \begin{equation}\label{counterexample}
        \| T^{\mathbf{d_r}} f \|_{L^2(\T)} \gtrsim \sqrt{\log \frac{1}{1-r}}  \| f\|_{L^2(\T)}\, .
    \end{equation}
\end{OldTheorem}
The upper bound \eqref{upperbound compact} of Theorem \ref{oldmain} is a perturbation/generalization of the classical Carleson theorem, as the MT system can be thought of as a perturbation of the trigonometric system.

The difference between Theorem \ref{main} and the lower bound \eqref{counterexample} is due to the logarithmic gap in the sequence $\mathbf{d_r}$ that is absent in $\mathbf{a_r}$. For \eqref{counterexample}, this logarithmic gap obstructs a nontrivial accumulation of phase and, hence, makes it impossible for any cancellation effects to emerge.
On the other hand, the operator $T^{\mathbf{a}_r}$ for the sequence $\mathbf{a}_r$ of Theorem \ref{main} is in some sense between Carleson type operators and singular integral operators. It accumulates enough phase to go beyond the singular integral theory and for cancellation effects of $TT^*$ type to occur and at the same time maintains good space localization not to fall into the world of the Carleson type theorems.

As pointed out in \cite{MNATSAKANYAN2022109461}, our problem is conformally invariant, i.e. for any M\"obius transform of the disk $m$, the boundedness of $T^{(a_j)_{j=1}^{\infty}}$ in $L^2$ is equivalent to the boundedness of $T^{(m(a_j) )_{j=1}^{\infty} }$. The M\"obius transform preserves hyperbolic distances. In the case of $\mathbf{d_r}$ the distance between any two points is $\sim \log \frac{1}{1-r}$. In contrast, the diameter of $\mathbf{a}_r$ is again $\log \frac{1}{1-r}$ but the distance between consecutive points is $\sim 1$. It seems to us that having fixed gaps between consecutive points should ensure boundedness of $T^{\mathbf{a}}$. We state this as a conjecture which generalizes Theorems \ref{main} and \ref{oldmain}.

\begin{conjecture}
Let $D>0$ be some constant, $\mathbf{a} = (a_n)_{n=1}^\infty$ be a sequence such that the hyperbolic distance between $a_n$ and $a_{n+1}$ is at most $D$ for any $n=1,2,\dots$. Then,
$$
\| T^{\mathbf{a}} \|_{L^2 (\T)\to L^2(\T)} \lesssim_D 1 \, .
$$
\end{conjecture}

The proof of Theorem \ref{main} is based on a $TT^*$ argument, the translational symmetry of the sequence $\mathbf{a}_r$ and the specific logarithmic accumulation of phase in the partial sum operator. The paper is organized as follows. In section \ref{preliminaries}, we recall and prove some preliminary lemmas. Section \ref{model case} presents a model case for Theorem \ref{main} illustrating and motivating the main idea of the proof. The proof of Theorem \ref{main} is presented in section \ref{proofofmain} and the proof of Corollary \ref{aeconvergence} in section \ref{proofofcor}. In section \ref{simple counterexample}, we present a simple proof of the lower bound \eqref{counterexample} of Theorem \ref{oldmain}. Unlike the proof in \cite{MNATSAKANYAN2022109461}, we obtain a pointwise lower bound for $T^{\mathbf{d_r}}f$.

\subsection*{Notation}
We will write  $A\lesssim B$ or $A=O(B)$ if $A\leq C B$ for some absolute constant $C>0$. We write $A\sim B$, if $A\lesssim B$ and $B\lesssim A$. For two functions $f,g$ defined on the unit circle $\T$, we denote by $\langle f,g \rangle := \frac{1}{2\pi} \int_{-\pi}^\pi f(e^{i\theta}) \overline{g(e^{i\theta})} d\theta$ their inner porduct in $L^2(\T)$. We will shorthand $\|f\|_2 := \| f\|_{L^2(\T)}$.


\subsection*{Acknowledgments}
This work was supported by the Higher Education and Science Committee of RA (Research Project No 24RL-1A028).

The argument presented in section \ref{simple counterexample} was suggested by Grigori Karagulyan whom I thank for many interesting discussions on orthogonal series.

\section{Preliminaries}\label{preliminaries}


Let $L\in \mathbb{N}$ and let $\mathbf{a} = (a_n)_{n=1}^{L}$ be an arbitrary finite sequence of points in the unit disk.

We start with notation. Let $\Psi_w$ be the phase of the M\"obius transform taking $w\in \D$ to $0$, i.e.
\begin{equation}
    \Psi_w (x) := x - \arg(w) + 2 \arcsin \frac{ |w| \sin ( x - \arg (w) ) }{ \sqrt{ 1 + |w|^2 - 2 |w| \cos( x - \arg (w))} } \, ,
\end{equation}
and
$$
\frac{\bar w}{|w|} \frac{e^{ix}-w}{1-\bar w e^{ix}} = e^{i\Psi_w (x)} \, .
$$
Denote also, for $1\leq n\leq L$,
\begin{equation}
    \psi_n^{\mathbf{a}} (x) = \psi_n (x) := \sum\limits_{j=1}^n \Psi_{a_j} (x).
\end{equation}

The Hardy-Littlewood maximal function on the circle will be denoted by $M$. Two versions of circular Hilbert transform are defined by
$$
Hf (x) = \int_\mathbb{T} f(y) \frac{dy}{\tan \frac{x-y}{2}} \, , \text{ and }\ti H f(x) = \int_\mathbb{T} f(y) \frac{dy}{\sin \frac{x-y}{2}} \, .
$$
Letting $\varphi (t) = \sin \frac{t}{2}$, we have
\begin{equation}\label{difference of hilbert transforms}
    \ti H f(x) - Hf(x) = 2 \int_\T f(y) \sin \frac{x-y}{2} dy = 2 (f * \varphi ) (x) \, .
\end{equation}
We know that for $f\in H^2(\D)$, $H f = - i (f - \int_\T f)$. We will use this property later when $f$ is a Blaschke product.

\begin{lemma}\label{partial sum esimate}
    Let $1\leq n \leq L$, then for the $n$th partial sums of the MT series,
    \begin{equation}\label{est max partial sum}
        \left| \sum_{j=0}^{n-1} \langle f, \phi_j \rangle \phi_j(e^{i\theta}) \right| \lesssim \left| \ti H \left(f e^{-i\psi_n(\cdot )} \right) (x) \right| + |H f(x)| + Mf(x) \, .
    \end{equation}
\end{lemma}

\begin{proof}
The finite dimensional subspace of $H^2(\D)$ spanned by $\phi_j$, $j=0,\dots, n-1$ is the orthogonal complement of the invariant subspace $B_nH^2(\D)$ \cite{Coif2}. Hence, the partial sum operator can be expressed as
$$\Id - \Proj_{B_nH^2}=\Id - B_n\Proj_{H^2}B_n^{-1}\, ,$$
where $\Id$ is the identity operator on $L^2(\T)$ and $\Proj_{V}$ is the orthogonal projection operator from $L^2(\T)$ on its subspace $V$. We have
\begin{equation*}
    \sum_{j=0}^{n-1} \langle f, \phi_j \rangle \phi_j(e^{i\theta}) = \int_\T f(y) \frac{B_n(e^{iy})^{-1}B_n(e^{ix})-1}{e^{i(x-y)}-1} dy \, .
\end{equation*}
$$
= \int_\T f(y) \frac{e^{i(\psi_n(x) - \psi_n(y))} - 1}{e^{i(x-y)}-1} dy
$$
$$
= e^{i\psi_n(x)} \int_\T f(y) \frac{e^{-i\psi_n(y)}dy}{e^{i(x-y)}-1} - \int_\T f(y) \frac{dy}{e^{i(x-y)}-1} \, .
$$
From the last expression \eqref{est max partial sum} follows by a triangle inequality and by changing the kernel $\frac{1}{e^{i(x-y)}-1}$ to $\frac{1}{\sin \frac{x-y}{2}}$ at an expense of a Hardy-Littlewood maximal function.
\end{proof}
For the analysis of the M\"obius phases we recall Lemmas 6.1 and 3.1 from \cite{MNATSAKANYAN2022109461}.
\begin{lemma}\label{decomp}
For any $0< r <1$ and $y\in [-\pi,\pi]$, we have
\begin{equation}\label{derratior}
\frac{(1-r)}{y^2+(1-r)^2} \lesssim \Psi_{r} ' (y) \lesssim \frac{(1-r)}{y^2+(1-r)^2} \, .
\end{equation}
\end{lemma}

\begin{lemma}\label{phaseasymptotics}
For any $0< r<1$ and $0\leq x\leq 1$ we have
\begin{equation}\label{singlephaseaction}
    | \Psi_{r} (x) - \pi + \frac{1}{1+\frac{x}{1-r}} | \lesssim \frac{1}{1+\Big(\frac{x}{1-r}\Big)^2} + \frac{1-r}{1+\Big(\frac{x}{1-r}\Big)}+(1-r) x.
\end{equation}
\end{lemma}
The estimate \eqref{singlephaseaction} is slightly better than the one in \cite[Lemma 6.1]{MNATSAKANYAN2022109461}, however, the proof is exactly the same.
Denote
$$
k(x) := \left[ \frac{x}{2\pi(1-r)} \right] \, .
$$
Then, \eqref{singlephaseaction} can equivallently be written in terms of $k(x)$, namely,
\begin{equation}
    \left| \Psi_{r} (x) - \pi + \frac{1}{2\pi(1+k(x))} \right| \lesssim \frac{1}{1+k(x)^2} + \frac{1-r}{1+k(x)}+(1-r)^2 k(x).
\end{equation}


Next, we want to introduce the linearized version of the Carleson operator \eqref{carleson op intro}. 
Let $N : \T \to \{1,2,\dots, L \}$ be an arbitrary measurable function. We put
$$
N(x) = \sum_{m=1}^L m \Z1_{A_m}(x) \, ,
$$
for pairwise disjoint measurable sets $A_m\subset \T$.
Then, we denote
\begin{equation}\label{linearized Carleson}
    T_N^{\mathbf{a}} f(x) = T_Nf(x) := \ti H \left( f e^{-i\psi_N(x)(\cdot )} \right) (x) = \int_\T f(y) \frac{e^{-i\psi_{N(x)} (y)}dy}{\sin \frac{x-y}{2}} \, .    
\end{equation}
We omit the superscript $\mathbf{a}$ in $T_N$ whenever it is clear from the context.

The last integral in \eqref{linearized Carleson} is well-defined for almost every $x$ as a finite sum of Hilbert transforms, that is,
$$
T_Nf(x) = \sum_{m=1}^L \ti H (f e^{-i\psi_m (\cdot)} ) \Z1_{A_m} (x) \, .
$$
Let us see what the adjoint $T_N^*$ of $T_N$ looks like. We have
$$
T_N^* g(y) = - \int_\T g(z) \frac{e^{i\psi_{N(z)}(y)}}{\sin \frac{z-y}{2}} dz = - \sum_{m=1}^L e^{i\psi_m(y)} \int_\T g(z)\Z1_{A_m}(z) \frac{dz}{\sin\frac{z-y}{2}}$$
$$
= - \sum_{m=1}^L e^{i\psi_m(y)} \ti H (g\Z1_{A_m}) (y) \, .
$$
The following lemma gives an integral formula for working with $T_NT_N^*$.
\begin{lemma}\label{TT* lemma}
    Let $g \in L^2(\T)$ be real valued. Then,
    \begin{equation}\label{TTstar formula}
        \|T^*_Ng\|_2^2  =O(\|g\|_2^2)+ 2 \int\limits_{\T ^2} g(x) g(z)  \frac{\sgn(N(z)-N(x)) \sin (\psi_{N(z)} - \psi_{N(x)} ) (x) }{\sin\frac{x-z}{2}} dxdz
    \end{equation}
    $$
    = O(\|g\|_2^2)+
    $$
\begin{equation}\label{TT*formula abs integrable}
         \int\limits_{\T^2} g(x) g(z)  \sgn(N(z)-N(x)) \frac{\sin (\psi_{N(z)} - \psi_{N(x)} ) (x) - \sin (\psi_{N(z)} - \psi_{N(x)} ) (z)}{\sin\frac{x-z}{2}} dxdz
    \end{equation}
\end{lemma}
The integrand function in \eqref{TT*formula abs integrable} is absolutely integrable. The integral on the right side of \eqref{TTstar formula} is obtained by \eqref{TT*formula abs integrable} by breaking the sine fraction and using the symmetry in $x$ and $z$. 
On the other hand, as $N$ attains finitely many values, the integral on the right side of \eqref{TTstar formula} can be represented as a finite sum of inner products of Hilbert transforms of $L^2$ functions with $L^2$ functions, hence, is well-defined.

The proof of the lemma relies on Fubini's theorem and the fact that the circular Hilbert transform $H$ acts on Blaschke products as $-i$ times the identity operator. However, as the Hilbert transforms of $L^2$ functions are defined for almost every point as a principal value integral, we will be careful in justifying the application of Fubini's theorem and this makes the proof of the lemma slightly long.

\begin{proof}
We have
$$
T_NT_N^* g (x) = -\sum_{m=1}^L\sum_{m'=1}^L \Z1_{A_m} (x) \ti H \left( e^{i(\psi_{m'} -\psi_m)(\cdot )} \ti H (g\Z1_{A_{m'}}) \right) (x)
$$
$$
=:I_1(x) + I_2(x) \, ,
$$
where $I_1$ corresponds to the diagonal sum of indices $m=m'$ and $I_2$ corresponds to the off-diagonal sum of indices $m\neq m'$.

For the first term we have
$$
|\langle I_1 ,  g \rangle| = \left| \langle \sum_{m=1}^L \Z1_{A_m} \ti H \left( \ti H (g\Z1_{A_{m}}) \right) , g \rangle \right|
= \sum_{m=1}^L \| \ti H (g\Z1_{A_m}) \|_2^2 \lesssim \| g\|_2^2.
$$

Let us turn to the second term.
We can write
$$
I_2(x) = -\sum_{m\neq m'} \Z1_{A_m} (x) \ti H \left( e^{i(\psi_{m'} -\psi_m)(\cdot )} \ti H (g\Z1_{A_{m'}}) \right) (x)
$$
$$
= -\sum_{m\neq m'} \Z1_{A_m} (x) e^{i(\psi_{m'} -\psi_m)(x)} \ti H \left(  \ti H (g\Z1_{A_{m'}}) \right) (x)
$$
$$
- \sum_{m\neq m'} \Z1_{A_m} (x) \ti H \left( (e^{i(\psi_{m'} -\psi_m)(\cdot )} - e^{i(\psi_{m'} -\psi_m)(x)}) \ti H (g\Z1_{A_{m'}}) \right) (x)
$$
$$
=: J_1(x) + J_2(x) \, .
$$

$J_2$ will be our main term. So we treat $J_1$ first. Let $K_h$ be the operator of convolution with the function $h$, i.e. $K_h f = f * h$. Then, by \eqref{difference of hilbert transforms},
$$
\ti H (\ti H (g\Z1_{A_{m'}})) = (H+2K_\varphi) (H +2K_\varphi) (g\Z1_{A_{m'}})
$$
$$
= H(H(g\Z1_{A_{m'}})) + 4 K_{\varphi *\varphi} (g\Z1_{A_{m'}}) + 4 H(\varphi) * (g\Z1_{A_{m'}})
$$
$$
= -g\Z1_{A_{m'}} +\int_\T g\Z1_{A_{m'}} + 4 K_{\varphi *\varphi} (g\Z1_{A_{m'}}) + 4 H(\varphi) * (g\Z1_{A_{m'}}) \, .
$$
Here, we used the fact $H(H(h)) = -h + \int_\T h$. Thus, using the fact that $A_m$ and $A_{m'}$ are disjoint, we have
$$
|\Z1_{A_m}\ti H (\ti H (g\Z1_{A_{m'}})) | \lesssim \Z1_{A_m} \left( \int_\T |g|\Z1_{A_m'} + |H(\varphi)| * (|g|\Z1_{A_{m'}})\right) \, .
$$
And we conclude our estiamte for $J_1$.
$$
|\langle J_1, g\rangle |
\lesssim \sum_{m\neq m'} \langle \Z1_{A_m}\left( \int_\T |g|\Z1_{A_m'} + |H(\varphi)| * (|g|\Z1_{A_{m'}})\right),|g| \rangle
$$
$$\leq \int_\T |g(x)| \Big( \int_\T |g| + |H(\varphi)|*|g| (x))  \Big)dx
$$
$$
\leq \Big(\int_T |g|\Big)^2 + \|g\|_2 \cdot \| |H(\varphi)|*|g| \|_2 \leq \|g\|_2^2 + \|g\|_2^2 \|H(\varphi)\|_1 \, .
$$
Here, we used Cauchy-Schwarz and Young's convolution inequality. As $\|H(\varphi)\|_1 \leq \| H(\varphi) \|_2^2 \lesssim \|\varphi\|_2^2\lesssim 1$, the estimate for $J_1$ is complete.

Finally, we turn to $J_2$. Let us denote the general term in $J_2$ by $J^{m,m'}_2$. We have
$$
-J^{m,m'}_2 (x) = \Z1_{A_m} (x) \ti H \left( (e^{i(\psi_{m'} -\psi_m)(\cdot )} - e^{i(\psi_{m'} -\psi_m)(x)}) \ti H (g\Z1_{A_{m'}}) \right) (x)
$$
\begin{equation}\label{second line limit delta}
    =\Z1_{A_m} (x) \int\limits_{\T} \left( \lim_{\delta \to 0} \int\limits_{|z-y|>\delta} g(z) \Z1_{A_{m'}}(z) \frac{dz}{\sin \frac{z-y}{2}} \right) \frac{e^{i(\psi_{m'} -\psi_m)(y)} - e^{i(\psi_{m'} -\psi_m)(x)}}{\sin \frac{x-y}{2}} dy \, .    
\end{equation}
The limit, as $\delta\to 0$, holds in $L^2$ sense, whereas the other factor is bounded, i.e.
$$\left| \frac{e^{i(\psi_{m'} -\psi_m)(y)} - e^{i(\psi_{m'} -\psi_m)(x)}}{\sin \frac{x-y}{2}} \right| \leq \sup_{\T} | (\psi_{m'} -\psi_m)'| \lesssim \frac{L}{\min_{1\leq n\leq L} (1-|a_n|)} \, ,
$$
where the last inequality follows from Lemma \ref{decomp}.
So, by the dominated convergence theorem, we continue
$$
\eqref{second line limit delta} =\Z1_{A_m} (x) \lim_{\delta\to 0} \int\limits_{\T} \left( \int\limits_{|z-y|>\delta} g(z) \Z1_{A_{m'}}(z) \frac{dz}{\sin \frac{z-y}{2}} \right) \cdot
$$
$$\frac{e^{i(\psi_{m'} -\psi_m)(y)} - e^{i(\psi_{m'} -\psi_m)(x)}}{\sin \frac{x-y}{2}} dy \, ,
$$
then by Fubini's theorem,
$$
= \Z1_{A_m} (x) \lim_{\delta\to 0} \int\limits_{\T} g(z) \Z1_{A_{m'}}(z) \left( \int\limits_{|z-y|>\delta}  \frac{e^{i(\psi_{m'} -\psi_m)(y)} - e^{i(\psi_{m'} -\psi_m)(x)}}{\sin \frac{x-y}{2} \sin \frac{z-y}{2}} dy \right) dz \, .
$$
We want to compute the big bracket above.
$$
A_{m,m'}^\delta (x,z) :=  \int\limits_{|z-y|>\delta}  \frac{e^{i(\psi_{m'} -\psi_m)(y)} - e^{i(\psi_{m'} -\psi_m)(x)}}{\sin \frac{x-y}{2} \sin \frac{z-y}{2}} dy =
$$
$$
\frac{1}{\sin \frac{x-z}{2}} \int\limits_{|z-y|>\delta} (e^{i(\psi_{m'} -\psi_m)(y)} - e^{i(\psi_{m'} -\psi_m)(x)}) \left( \frac{1}{\tan \frac{x-y}{2}} - \frac{1}{\tan \frac{z-y}{2}} \right) dy
$$
$$
= \frac{1}{\sin \frac{x-z}{2}} \left( -H (e^{i(\psi_{m'} -\psi_m)}) (z) - \int_{|z-y|<\delta} \frac{e^{i(\psi_{m'} -\psi_m)(y)}-e^{i(\psi_{m'} -\psi_m)(z)}}{\tan\frac{z-y}{2}} dy\right)
$$
$$
+\frac{1}{\sin \frac{x-z}{2}} \left( H (e^{i(\psi_{m'} -\psi_m)}) (x) - \int_{|z-y|<\delta} \frac{e^{i(\psi_{m'} -\psi_m)(y)}-e^{i(\psi_{m'} -\psi_m)(x)}}{\tan\frac{x-y}{2}} dy\right) \, .
$$
We know that $e^{i(\psi_m-\psi_{m'})}$ is either a Blaschke product or its conjugate. Hence, 
$$H(e^{i(\psi_m-\psi_{m'})}) (x) = -i \sgn(m-m') \Big( e^{i(\psi_m-\psi_{m'})(x)} - \int_\T e^{i(\psi_m-\psi_{m'})} \Big) \, . $$
Thus, we can estimate
$$
\left| A^\delta_{m,m'} (x,z) - \frac{i\sgn (m-m')}{\sin\frac{x-z}{2}} \Big(e^{i(\psi_m-\psi_{m'})(z)} - e^{i(\psi_m-\psi_{m'})(x)}  \Big) \right|
$$
$$
\lesssim \frac{\delta}{|\sin \frac{x-z}{2}|} \cdot \frac{L}{\min_{1\leq n \leq L} (1- |a_n|)} \, .
$$
If $A_m$ and $A_{m'}$ are  unions of finitely many open intervals and $x\in A_m$, then
$$
\int_\T |g(z)| \Z1_{A_{m'}} (z) \frac{dz}{|\sin \frac{z-x}{2}|} < \infty
$$
for all but possibly finitely many $x\in A_m$ that are the common endpoints of the open intervals of $A_m$ and $A_m'$. Thus, for all but finitely many $x\in A_m$, by the dominated convergence theorem, we have
$$
-J_2^{m,m'}(x) = \Z1_{A_m} (x) i\sgn (m-m') \int\limits_{\T} g(z) \Z1_{A_{m'}}(z) \frac{e^{i(\psi_m-\psi_{m'})(z)} - e^{i(\psi_m-\psi_{m'})(x)}}{\sin\frac{x-z}{2}} dz \, .
$$
We can rewrite this as
\begin{equation}\label{main identity in this lemma}
    -J_2^{m,m'} (x) = \Z1_{A_m} (x) i\sgn (m-m') \ti H\left( g\Z1_{A_{m'}} (e^{i(\psi_m-\psi_{m'})(\cdot)} - e^{i(\psi_m-\psi_{m'})(x)}) \right) \, .
\end{equation}
To see that \eqref{main identity in this lemma} holds for general measurable sets, one just needs to approximate $A_m$ and $A_m'$ with finite unions of disjoint open intervals and pass to an $L^2$ limit on both sides of the equality. As $J_2^{m,m'}$ consists of two Hilbert transforms, the $L^2$ limits on both sides exist.

Recalling that $g$ is real valued and using \eqref{main identity in this lemma}, we compute
$$
\langle J_2, g\rangle = \sum_{m\neq m'} \langle g, J_2^{m,m'} \rangle
$$
$$
=-\sum_{m\neq m'} i\sgn (m-m') \int\limits_\T\int\limits_\T g(x)g(z)\Z1_{A_m} (x) \Z1_{A_{m'}} (z) \cdot
$$
$$\frac{e^{i(\psi_m-\psi_{m'})(z)} - e^{i(\psi_m-\psi_{m'})(x)}}{\sin\frac{x-z}{2}} dzdx \, .
$$
Using the symmetry between $m$ and $m'$, we continue.
$$
=-\frac{i}{2} \sum_{m\neq m'} i\sgn (m-m') \int\limits_\T\int\limits_\T g(x)g(z)\Z1_{A_m} (x) \Z1_{A_{m'}} (z) \cdot
$$
$$\frac{e^{i(\psi_m-\psi_{m'})(z)} - e^{i(\psi_m-\psi_{m'})(x)}}{\sin\frac{x-z}{2}} dzdx
$$
$$
-\frac{i}{2} \sum_{m\neq m'} (-\sgn (m-m') ) \int\limits_\T\int\limits_\T g(x)g(z)\Z1_{A_{m'}} (x) \Z1_{A_{m}} (z) \cdot
$$
$$\frac{e^{-i(\psi_m-\psi_{m'})(z)} - e^{-i(\psi_m-\psi_{m'})(x)}}{\sin\frac{x-z}{2}} dzdx \, ,
$$
then switching $x$ and $z$ in the second double integral we get,
$$
= \sum_{m\neq m'} i\sgn (m-m') \int\limits_\T\int\limits_\T g(x)g(z)\Z1_{A_m} (x) \Z1_{A_{m'}} (z) \cdot
$$
$$\frac{\sin ((\psi_m-\psi_{m'})(z)) - \sin ((\psi_m-\psi_{m'})(x))}{\sin\frac{x-z}{2}} dzdx \, .
$$
Finally, separating the two sine terms and closing the summation in $m, m'$, we get
$$
= 2 \int_\T \int_\T g(x) g(z) \sgn(N(z)-N(x)) \sin \left( (\psi_{N(z)} - \psi_{N(x)} ) (x) \right) \frac{ dxdz}{\sin\frac{x-z}{2}} \, .
$$

\end{proof}

Let us make some notation for the modulation factor and the quadratic form in \eqref{TTstar formula}. Let
\begin{equation}\label{modulation chi}
    \chi_N (x,z) := \sgn (N(z)-N(x)) \sin ( (\psi_{N(z)} - \psi_{N(x)}) (x) ) \, ,
\end{equation}
and
\begin{equation}\label{quadraticform B}
    B(g, N) := \int_\T \int_\T g(x) g(z)  \frac{\chi_N (x,z)}{\sin\frac{z-x}{2}} dxdz \,.
\end{equation}
As $\langle T_NT_N^* g,g\rangle\geq 0$, the lemma implies, that
\begin{equation}\label{positive quadratic form}
    B(g,N) \gtrsim - \|g\|_2^2 \, .
\end{equation}

\section{A model case for the proof of Theorem \ref{main}}\label{model case}
Let $I_m = [2\pi (1-r), 2\pi (1-r)(m+1)]$ for $m\in \ZZ$.
By Lemma \ref{TT* lemma}, Theorem \ref{main} follows if we can show $B(g,N)\leq \|g\|_2^2$. We want to give a heuristic proof of the latter for the following $g$ and $N$. Let $g(x) = \alpha_m$ for $x\in I_m$, $m$ is even and $0\leq m\leq [1/(1-r)]$, and $\alpha_m$'s are arbitrary real numbers, and is $0$ elsewhere. $N(x) = k(x)$.

The heuristic argument below demonstrates the use of $TT^*$ and the symmetries and motivates the main construction in subsection \ref{aux constructions}. Let us plug $g$ and $N$ into \eqref{quadraticform B}. We have
$$
B(g, N) = \sum_{\substack{0\leq j,j'\leq [1/(1-r)] \\ j\neq j', \text{ even} }} \alpha_j\alpha_{j'} \int_{I_j}\int_{I_{j'}} \frac{\sgn (j'-j) \sin (\psi_{j'}-\psi_j)(x) }{\sin{\frac{z-x}{2}}} dxdz \, .
$$
We can approximate the last expression as follows. Let $x\in I_j, z\in I_{j'}$ and $j<j'$, then
$$\sin\frac{z-x}{2} \approx \frac{z-x}{2} \approx \pi(1-r)(j'-j) \, .$$
Next, by Lemma \ref{phaseasymptotics},
$$
(\psi_{j'}-\psi_j) (x) \approx \sum_{k=j+1}^{j'} \pi - \frac{1}{2\pi (k-j)} = (j'-j)\pi - \sum_{k=1}^{j'-j} \frac{1}{2\pi k} \approx (j'-j)\pi - \frac{\log (j'-j)}{2\pi} \, .
$$
Combining the two, we write
$$
B(g,N) \approx - 4\pi \sum_{\substack{0\leq j <j'\leq [1/(1-r)] \\ j,j' \text{ even} }} \alpha_j\alpha_{j'} \frac{\sin (\frac{1}{2\pi} \log (j'-j) )}{j'-j} (1-r) \, .
$$
We have $\|g\|_{2}^2 = \sum\limits_{\substack{0\leq j\leq [1/(1-r)] \\ j \text{ even}}} \alpha_j^2 (1-r)$, so we need to show
$$
-\sum_{\substack{0\leq j <j'\leq [1/(1-r)] \\ j,j' \text{ even} }} \alpha_j\alpha_{j'} \frac{\sin (\frac{1}{2\pi} \log (j'-j) )}{j'-j} \lesssim \sum\limits_{\substack{0\leq j\leq [1/(1-r)] \\ j \text{ even}}} \alpha_j^2 \, .
$$

For an arbitrary real sequence $\alpha = (\alpha_j )_{j=0}^\infty$, let
$$
T (\alpha) := - \sum_{0\leq j < j'} \alpha_j\alpha_{j'} \frac{\sin (\frac{1}{2\pi} \log (j'-j) )}{j'-j} \, .
$$
The inequality \eqref{positive quadratic form} means that we can assume $T(\alpha)\gtrsim -\|\alpha\|_{\ell^2}^2$. We want to show that $T(\alpha) \lesssim \|\alpha\|_{\ell^2}^2$. We construct a sequence $\beta$ in the following way $\beta_{[e^{2\pi^2 } j] } = \alpha_j$ and $\beta_k = 0$ if $k\neq [e^{2\pi^2} j]$ for any $j\in \ZZ$. Then, $\|\beta \|_{\ell^2 } = \|\alpha\|_{\ell^2}$, and we can compute
$$
T(\beta) = \sum_{0\leq j < j'} \alpha_j\alpha_{j'} \frac{\sin (\frac{1}{2\pi} \log ([e^{2\pi^2}j']-[e^{2\pi^2}j ]) ) }{[e^{2\pi^2}j']-[e^{2\pi^2}j ]}
$$
$$
\approx \sum_{0\leq j < j'} \alpha_j\alpha_{j'} \frac{\sin (\frac{1}{2\pi} \log (e^{2\pi^2} (j'-j)) ) }{e^{2\pi^2}(j'-j)}
$$
$$
=\sum_{0\leq j < j'\leq R} \alpha_j\alpha_{j'} \frac{\sin (\pi + \frac{1}{2\pi} \log (j'-j) )}{e^{2\pi^2}(j'-j)} \approx - \frac{1}{e^{2\pi^2}} T(\alpha) \, .
$$
As $T(\beta) \gtrsim -\|\alpha\|_{\ell^2}^2$ and by the last computation $T(\beta) \approx - \frac{1}{e^{2\pi^2}} T(\alpha)$, we conclude $T(\alpha) \lesssim \|\alpha\|_{\ell^2}^2$ as wanted.s

In the presentation above $N$ was chosen so as to maximize the effect of the Blaschke modulation factors at each point. In general, when $N(x)$ is not necessarily close to $k(x)$ we will need an additional construction reflecting $N(x)$ across $k(x)$. Then, in the proof of the general case presented in the following section, for a given pair $(g,N)$, we construct 3 additional pairs. The first one is similar to the construction of $\beta$ above, the second one is a reflection of $(g,N)$ to compensate for the fact that $N(x)$ may not be close to $k(x)$, and the last one is "the $\beta$ construction" for the reflected pair.

\section{Proof of Theorem \ref{main}}\label{proofofmain}
This section is organized as follows. In the first part we reduce the proof of Theorem \ref{main} to Claim \ref{claim1} and Claim \ref{claim2}. In subsection \ref{aux constructions}, we provide some auxiliary constructions and lemmas. Then, Claim \ref{claim1} is proved in subsection \ref{proof of claim1} and Claim \ref{claim2} is proved in subsection \ref{proof of claim2}.

Let $K = \left[\frac{1}{4(1-r)} \right]$. We extend the sequence $\mathbf{a_r}$ to have also negative indices. That is, we put $a_n=re^{2\pi i n(1-r)}$, for $-K\leq n\leq 0$. Then, we can consider the MT system corresponding to the new sequence $\mathbf{a_r} = (a_n)_{n=-K}^{[1/(1-r)]}$. Namely, let $\phi_{-K-1}(z) = \frac{\sqrt{1-|a_{-K}|^2}}{1-\overline{a_{-K}}z}$ and, for $-K\leq n \leq \left[\frac{1}{1-r} \right]-1$,
$$
B_n(z) = \prod_{j=-K}^n \frac{\bar a_j}{|a_j|}\frac{z-a_j}{1-\overline{a_j} z} \text{ and } \phi_{n} (z) = B_n(z) \frac{\sqrt{1-|a_{n+1}|^2}}{1- \overline{a_{n+1}}z} \, .
$$
As before, we denote
$$\psi_m (x) := \sum_{j=-K}^m \Psi_{a_j}(x) \, , -K\leq m \leq K \, .$$
Next, let
$$N: \T \to  \left\{ -\left[\frac{e^{-2\pi^2}}{16(1-r)}\right],\dots, 0,1,\dots, \left[\frac{e^{-2\pi^2}}{16(1-r)}\right] \right\}$$
be an arbitrary measurable function. Recall the definition \eqref{linearized Carleson},
$$
T_N f(x) = \int_\T f(y) e^{-i\psi_{N(x)}(y)} \frac{dy}{\sin\frac{x-y}{2}} \, .
$$
We denote
$$p_N(x)=p(x):= N(x)-k(x) \, .$$
And let
\begin{equation}
    E_N = E := \left\{x\in \left[-e^{-2\pi^2}/2, e^{-2\pi^2}/2 \right] \, : \, k(x) \text{ is even and } p(x)\geq 0 \right\} \,.
\end{equation}

Then, our first claim is as follows.
\begin{claim}\label{claim1}
    If for any $N$ as above and any real-valued $g\in L^2(\T)$ supported on $E$, we have
    \begin{equation}\label{claim1 estimate}
        \| T_N^* g\|_2 \lesssim \|g\|_2 \, ,
    \end{equation}
    then Theorem \ref{main} holds.
\end{claim}
The proof of Claim \ref{claim1} is rather straightforward if technically somewhat unpleasant.

In subsection \ref{aux constructions}, we construct, for any pair of functions $(g,N)$ satisfying the hypothesis of Claim \ref{claim1}, three pairs of functions $(f,M), (\ti g, \ti N), (\ti f, \ti M)$, such that
$$M, \tilde{M}, \tilde{N} \,: \, \T \to \{-K, -K +1,\dots, K-1, K \}$$
are measurable, and
$$f,\tilde{f}, \tilde{g} \,:\, \T \to \R \, $$
with
$$
\|g\|_2=\|f\|_2=\|\ti g\|_2=\|\ti f\|_2 \, .
$$
Then, our second claim is as follows.
\begin{claim}\label{claim2}
    For an arbitrary pair of functions $(g,N)$ satisfying the hypothesis of Claim \ref{claim1}, and the three associated pairs $(\ti g, \ti N), (f, M), (\ti f, \ti M)$, we have
    \begin{equation}\label{claim2 estimate}
        B(g,N) + e^{2\pi^2} B(\ti g,\ti N) + B(f,M) + e^{2\pi^2} B(\ti f, \ti M) \lesssim \|g\|_2^2 \, .
    \end{equation}
\end{claim}
As $L^2$ norms of $g,\ti g, f,\ti f$ are all equal, by the estimate \eqref{positive quadratic form}, Claim \ref{claim2} will imply
$$
B(g,N) \leq O(\|g\|_2^2) + B(g,N) + e^{2\pi^2} B(\ti g,\ti N) + B(f,M) + e^{2\pi^2} B(\ti f, \ti M) \lesssim \|g\|_2^2 \, .
$$
With Claim \ref{claim1}, this will prove Theorem \ref{main}.

\subsection{Constructing the auxiliary functions}\label{aux constructions}
Recall $I_m = [2\pi (1-r)m, 2\pi (1-r)(m+1)]$ for $m\in \mathbf{Z}$. Denote
$$
\tau (x) : = 2\pi(1-r)(2k(x)+1)-x \, .
$$
$\tau(x)$ is the reflection of $x$ with respect to the midpoint of $I_{k(x)}$. 

Denote
\begin{equation}
    \ti k (x) := \begin{cases}
    & [e^{2\pi ^2} k(x)], \text{ if it is even,} \\
    & [e^{2\pi ^2} k(x)] + 1, \text{otherwise}
    \end{cases} \, ,
\end{equation}
and
$$
\eta (x) := \left( \ti k(x) - k(x)\right)2\pi (1-r) + x \, .
$$
$\eta(x)$ is a sort of dilation and translation of $x$.
The following lemma states some properties of $\eta$ and $\tau$.
\begin{lemma}\label{taueta} For any $x,z\in [-\pi,\pi]$, we have

    \begin{itemize}
        \item[a.] $\tau(\eta(x))=\eta(\tau(x))$,
        \item[b.] $\tau^2(x) = x$ and $\tau^{-1}(x)= \tau(x)$,
        \item[c.] $k(\eta(x)) = \ti k(x)$,
        \item[d.] $\eta (x) - \eta (z) = e^{2\pi ^2}(x-z) + O(1-r)$,
        \item[e.] $\tau(x)-\tau(z) = x-z + O(1-r)$,
        \item[f.] $\eta$ is injective and $\eta ([-e^{-2\pi^2}/2,e^{-2\pi^2}/2]) \subset [-1,1]$.
    \end{itemize}
\end{lemma}

\begin{proof}
The proof is straightforward. Let us only check part f. For all $|x|\leq e^{-2\pi^2}/2$, we have
$$ |\eta (x)| \leq 2\pi(1-r) |\ti k(x)| + \left| x - 2\pi (1-r)k(x) \right| $$
$$\leq 2\pi(1-r)e^{2\pi ^2} \frac{|x|}{2\pi(1-r)} + 2\pi(1-r) \leq \frac{1}{2}+\frac{1}{2} \leq 1 \, . $$
This proves the inclusion. For injectivity, we have $\eta(x) \in I_{\ti k(x)}$. Hence, if $\eta(x) = \eta (x')$, then $\ti k(x) = \ti k(x')$. $k(x)=k(x')$ by definition of $\ti k$, and $x=x'$ by definition of $\eta$.
\end{proof}
The above lemma will be repeatedly used without a reference.

Let us fix the pair $(g,N)$ and construct the associated three pairs. Morally, $(M,f)$ is a reflection of $(N,g)$. $(\ti N,\ti g)$ is a dilation and translation of $(N,g)$, and $(\ti M , \ti f)$ is obtained from $(M,f)$ the same way as $(\ti N, \ti g)$ is obtained from $(N,g)$. We put
$$M(x) := \left( k(x) + 1 -p_N( \tau (x) ) \right) \Z1_{E}(\tau(x))$$
and
\begin{equation}
 f(x) := g \Big( \tau (x) \Big) \Z1_{E}( \tau(x)).
\end{equation}
Obviously $M$ is measurable and $\|f\|_2=\|g\|_2$. We need to check the range of $M$. For $x\in \tau(E)\subset [-1/48,1/48]$,
$$|M(x)| \leq |k(x)|+p_N (\tau(x) ) +1 \leq \left[ \frac{e^{-2\pi^2}}{4\pi(1-r)}\right] + 2\left[ \frac{e^{-2\pi^2}}{16(1-r)}\right] \leq K \, .$$
Next,
$$\ti N (x) := \left( \ti k( \eta^{-1}(x) ) + p_N(\eta^{-1}(x) ) \right) \Z1_{\eta(E)}(x) \, ,
$$
or in other words, for $x\in E$,
\begin{equation}
    \ti N \left( \eta (x) \right) = \ti k(x) + p_N(x) = \ti k(x) + N(x) - k(x) \, .
\end{equation}
And
$$\ti g (x) := g(\eta^{-1}(x) ) \Z1_{\eta(E)} (x) \, .
$$
Again, that $\|\ti g\|_2 = \|g\|_2$ and $\ti N$ is measurable is trivial. Let us check the range of $\ti N$. For any $x\in E$,
$$
|\ti N(\eta(x))| \leq |\ti k(x)| + |p_N(x)| \leq e^{2\pi^2} \left[ \frac{e^{-2\pi^2}/2}{2\pi \cdot (1-r)} \right] + 3\left[\frac{e^{-2\pi^2}}{16(1-r)}\right] \leq K \, .
$$

Then, the pair $(\tilde{M}, \tilde{f})$ is obtained from $(M,f)$ by the same procedure as $(\tilde{N}, \tilde{g})$ from $(N,g)$. Namely,
$$\ti M (x) := \left(\ti k(\eta^{-1} (x) ) + M(\eta^{-1}(x) ) - k (\eta^{-1}(x)) \right) \Z1_{\eta (\tau(E))} (x)
 \, ,
$$
and
$$\ti f (x) := f(\eta^{-1}(x)) \Z1_{\eta(\tau(E))} (x) \, .
$$
One can check that
$$\ti M (x) = \Big( k(x) + 1 - p_{\ti N} (\tau (x)) \Big) \Z1_{\eta (\tau(E) )} (x)
 \, ,
$$
and
$$\ti f (x) := g\Big( \eta^{-1}(\tau(x)) \Big) \Z1_{\eta(\tau(E))} (x) \, .
$$
So that $(\ti M, \ti f)$ can equivalently be obtained from $(\ti N,\ti g)$ the same way $(M,f)$ is obtained from $(N,g)$.

Let us finish this subsection with a lemma about the phases of the Blaschke products.
\begin{lemma}\label{psi phase estimate}
    For $-K\leq m<m'\leq K$ and $y\in [-\pi,\pi]$, we have
    $$
    0< (\psi_{m'} -\psi_{m})'(y) \lesssim \frac{1}{(1-r)(1+\dist (k(y), [m, m']))} \, . 
    $$
\end{lemma}
\begin{proof}
    By Lemma \ref{decomp},
    $$
    (\psi_{m'}-\psi_{m})'(y) = \sum_{j=m+1}^{m'} \Psi_{a_j}' (y)
    = \sum_{j=m+1}^{m'} \Psi_{r}' (y-\arg a_j) 
    $$
    $$= \sum_{j=m}^{m'} \Psi_{r}' (y-2\pi(1-r)j) \lesssim \sum_{j=m}^{m'} \frac{1-r}{(1-r)^2+(y-2\pi j(1-r))^2}
    $$
    $$
    \lesssim \frac{1}{(1-r)(1+\dist (\frac{y}{1-r}, [2\pi m,2\pi m']) )} \, .
    $$
\end{proof}

\subsection{Proof of Claim \ref{claim1}}\label{proof of claim1}
First, we want to drop the restrictions on the range and the support of $g$. By a triangle inequality, we can allow $g$ to be complex valued.

By rotational symmetry of the sequence $\mathbf{a_r}$, we can drop the evenness of $k(x)$. To see this, assume $g$ is an arbitrary function supported on
$$\{x\in \big[-e^{-2\pi^2}/2,e^{-2\pi^2}/2 \big] \, : \, p(x)\geq 0 \} \, .$$
We write
$$
g(x) = g(x)\Z1_{\{k(x) \text{ even}\}} +g(x) \Z1_{\{k(x) \text{ odd} \}} =: g_1(x)+g_2(x) \, .
$$
Then, \eqref{claim1 estimate} holds for $(g_1,N)$. On the other hand,
$$
T_N^*(g_2) (y) = - \int_\T g_2(z) \frac{e^{i\psi_{N(z)}(y)}}{\sin \frac{z-y}{2}} dz
$$
$$
= - \int_\T g_2(z-2\pi (1-r)) \frac{e^{i\psi_{N(z-2\pi (1-r))}(y)}}{\sin \frac{z-(y+2\pi (1-r))}{2}} dz
$$
$$
= - \int_\T g_2(z-2\pi (1-r)) \frac{e^{i\psi_{N(z-2\pi (1-r))}(y+2\pi (1-r) - 2\pi (1-r))} }{\sin \frac{z-(y+2\pi (1-r))}{2}} dz \, .
$$
By the rotational symmetry of $\mathbf{a_r}$, we have, for any $-K\leq m\leq K$,
$$
\psi_{m}(y+2\pi (1-r) - 2\pi (1-r)) = \Psi_{re^{-2\pi i (K+1)(1-r)}} (y+2\pi (1-r)) +  \psi_{m- 1}(y+2\pi (1-r)) \, ,
$$
Hence,
\begin{equation}
   |T_N^*(g_2) (y-2\pi (1-r))| =  \left| \int_\T g_2(z-2\pi (1-r))  \frac{e^{i\psi_{N(z-2\pi (1-r)) - 1}(y)} }{\sin \frac{z-y}{2}} dz \right| \, .
\end{equation}
The pair
$$\Big( g_2(\cdot - 2\pi (1-r)), N (\cdot - 2\pi (1-r)) -1 \Big)$$
satisfies the hypothesis of Claim \ref{claim1}, therefore,
$$
\|T_N^*(g_2)\|_2 \lesssim \|g_2\|_2 \, .
$$
Hence, we can apply a triangle inequality, to get
$$
\|T_N^* g\|_2\leq \|T_N^* g_1\|_2 + \|T_N^* g_2\|_2 \lesssim \|g\|_2 \, . 
$$

The condition $p(x)\geq 0$ is equivalent to $N(x)\geq k(x)$. To get rid of that, we use the reflection and modulation symmetries. Namely, let $g$ be a function with support
$$
\{x\in \big[-e^{-2\pi^2}/2,e^{-2\pi^2}/2 \big] \, : \, N(x) < k(x) \} \, .$$
For the Blaschke factors, we have
$$
e^{-i\Psi_w (-y)} = \frac{w }{|w|} \overline{\left(\frac{e^{-iy} - w}{1-\bar w e^{-iy}} \right) } = e^{i \Psi_{\overline{w}} (y)} \, .
$$
Hence, by the reflection symmetry of $\mathbf{a_r}$,
$$
e^{i\psi_m (-y)} = e^{-i\sum_{j=-K}^m \Psi_{\bar a_j} (y)} = e^{-i\psi_K (y)} \cdot e^{i\psi_{-(m+1)} (y)}$$
Then, we have
$$
T_N^* (g)(-y) = - \int_\T g(z) \frac{e^{i\psi_{N(z)}(-y)}}{\sin \frac{z+y}{2}} dz
= e^{-i\psi_K (y)} \int_\T g(-z) \frac{e^{i\psi_{-(N(-z)+1)}(y)}}{\sin \frac{z-y}{2}} dz \, .$$
We see that $N(-z)<k(-z)\leq  -k(z)-1$ implies
$-(N(-z)+1) > k(z)$. Hence, for the pair
$$\big( g(-\cdot), -(N(-\cdot) +1) \big)$$
the inequality \eqref{claim1 estimate} holds, and by another triangle inequality, we drop the condition $p(x)\geq 0$.

Finally, let us drop the last assumption on the support of $g$.
Assume $g$ is supported on $\T$. As
$$\dist \left( [-\pi,\pi]\setminus [-e^{-2\pi^2}/2,e^{-2\pi^2}/2], [-2\pi e^{-2\pi^2}/16, 2\pi e^{-2\pi^2}/2] \right) \sim 1 \, ,
$$
we have by Lemma \ref{psi phase estimate}, for all $-K\leq m, m'\leq K$ and $y\in [-\pi,\pi]\setminus [-e^{-2\pi^2}/2,e^{-2\pi^2}/2]$,
$$
|(\psi_m -\psi_{m'})'(y) | \lesssim 1 \, .
$$
Hence, using a triangle inequality, the inequality \eqref{claim1 estimate} for $g\Z1_{[-e^{-2\pi^2}/2,e^{-2\pi^2}/2]}$, and the mean value theorem applied to the formula \eqref{TT*formula abs integrable} of Lemma \ref{TT* lemma}, we have
$$
\| T^*_N (g ) \|_2 \lesssim \| g\|_2^2 +
\Big|\int\limits_{[-\pi , \pi]^2\setminus [-\frac{e^{-2\pi^2}}{2}, \frac{e^{-2\pi^2}}{2}]^2} g(x)g(z) \sgn (N(z)-N(x)) \cdot $$
$$\frac{\sin\big((\psi_{N(z)}-\psi_{N(x)})(x)\big) - \sin\big((\psi_{N(z)}-\psi_{N(x)})(z)\big) }{\sin\frac{x-z}{2}} dxdz \Big|
$$
$$
\lesssim \|g\|_2^2 + \int\limits_{[-\pi , \pi]^2\setminus [-\frac{e^{-2\pi^2}}{2}, \frac{e^{-2\pi^2}}{2}]^2} |g(x)g(z)| dxdz \lesssim \| g\|_2^2 \, .
$$

So far, we have shown that if \eqref{claim1 estimate} holds for any $(g,N)$ satisfying the hypothesis of Claim \ref{claim1}, then it also holds for any $N$ satisfying the same hypothesis and any $g\in L^2(\T)$. Next, let us make some reductions from Theorem \ref{main}.

Denote
$$\ti \phi_n := \phi_n^{\mathbf{a_r}|_{\left[1,[\frac{1}{1-r}] \right]}}\, ,\quad 0\leq n\leq \left[\frac{1}{1-r}\right] \, ,$$
the MT system corresponding to the restricted sequence $\mathbf{a_r}|_{[1,[\frac{1}{1-r}]]}=(a_n)_{n=1}^{[\frac{1}{1-r}]}$. Note that $\ti\phi_n$ is the original MT sequence that appears in the statement of Theorem \ref{main}. We have
$$
B_0(z) \ti \phi_j(z) = \phi_j(z)\, , \text{ for } 0\leq j \leq [1/(1-r)] \, .
$$
So that
$$
\sum_{j=1}^n \langle f,\ti\phi_j\rangle \ti \phi_j = \overline{B_0} \sum_{j=1}^n \langle fB_0 , \phi_j \rangle \phi_j \, .
$$
Hence, inequality \eqref{maininequality} is equivalent to
$$
\| \sup_{0\leq n \leq [1/(1-r)]} | \sum_{j=1}^n \langle f, \phi_j \rangle \phi_j | \|_2 \lesssim \|f\|_2 \, .
$$
By a triangle inequality
$$
\| \sup_{0\leq n \leq [1/(1-r)]} |\sum_{j=1}^n \langle f,\phi_j\rangle \phi_j | \|_2 \leq \sum_{l=0}^{16e^{2\pi^2}} \|\sup_{\frac{l}{16e^{2\pi^2}} \leq n \leq \frac{l+1}{16e^{2\pi^2}}} |\sum_{j=\frac{l}{16e^{2\pi^2}}}^n \langle f,\phi_j\rangle \phi_j | \|_2 \, .
$$
By the rotation symmetry of our sequence, it is enough to prove the boundedness only for the term corresponding to $l=0$ in the sum above.
As the supremums in the above displays are over finite sets, we can choose a function
$$N: \T \to  \left\{ 0,1,\dots, \left[\frac{e^{-2\pi^2}}{16(1-r)}\right] \right\} \, ,$$
depending on $f$, such that
$$
\sup_{0\leq n \leq [\frac{e^{-2\pi^2}}{16(1-r)} ]} |\sum_{j=0}^n \langle f,\phi_j\rangle \phi_j | = |\sum_{j=0}^{N(x)} \langle f,\phi_j\rangle \phi_j | \, .
$$
It is not difficult to see that $N$ is, in fact, a measurable function.
Therefore, Theorem \ref{main} is reduced to showing that for an arbitrary fixed measurable function $N: \T \to  \left\{ -\left[\frac{e^{-2\pi^2}}{16(1-r)}\right],\dots, \left[\frac{e^{-2\pi^2}}{16(1-r)}\right] \right\}$, and for all $f$, the inequality
\begin{equation}\label{one more reduction}
    \|\sum_{j=0}^{N(x)} \langle f,\phi_j\rangle \phi_j \| \lesssim \|f\|_2\, ,
\end{equation}
holds with the implicit constant independent of $N$. Then, Lemma \ref{partial sum esimate} finishes the proof of claim \ref{claim1}.


\subsection{The proof of Claim \ref{claim2}}\label{proof of claim2}
We want to prove \eqref{claim2 estimate}. Let us denote the left-hand side of \eqref{claim2 estimate} by $\Sigma$.

As the supports of $g,f,\ti g,\ti f$  are all inside $[-1,1]$, in the integral formula \eqref{quadraticform B} we can change the kernel from $\frac{1}{\sin\frac{x-z}{2}}$ to $\frac{1}{x-z}$ at an expense of a Hardy-Littlewood maximal function.
Then, by a change of variables, we write
$$
B(f, M) = O(\|g\|_2^2) + \int\limits_{\tau(E)\times \tau(E)} f(x) f(z) \chi_{M} (x,z) \frac{dxdz}{x-z}
$$
$$
= O(\|g\|_2^2) + \int\limits_{E\times E} g(x) g(z) \chi_{M} (\tau(x),\tau(z)) \frac{dxdz}{\tau(x)-\tau(z)} \, .
$$
Analogous formulas can be written for $B(\ti g, \ti N)$ and $B(\ti f, \ti M)$. Hence, we have
$$
\Sigma = O(\|g\|_2^2) + \int\limits_{E\times E} g(x)g(z) \Big( \frac{\chi_N (x,z)}{x-z} +\frac{e^{2 \pi ^2} \chi_{\tilde{N}} (\eta(x),\eta(z))}{\eta (x) -\eta (z)} + $$
\begin{equation}\label{Sigma kernel}
    \frac{\chi_M (\tau(x),\tau(z))}{\tau(x)-\tau(z)} + \frac{e^{2 \pi ^2} \chi_{\tilde{M}} (\eta(\tau(x)),\eta(\tau(z)) )}{\eta(\tau(x))-\eta(\tau(z))} \Big) dxdz  .  
\end{equation}

Let
$$
E_0 := \{ (x,z)\in E\times E \, : \, k(x)=k(z) \}
$$
be the diagonal. 
Then, by Lemma \ref{psi phase estimate}, we have
$$
\left| \int\limits_{E_0} g(x)g(z) \frac{\chi_N (x,z)}{x-z} dxdz \right| = \frac{1}{2} \Big| \int\limits_{E_0} g(x)g(z) \sgn (N(x)-N(z)) \cdot
$$
$$\frac{\sin \big( (\psi_{N(x)}-\psi_{N(z)})(z) \big) - \sin \big( (\psi_{N(x)}-\psi_{N(z)})(x) \big) }{x-z} dxdz \Big| $$
$$
\lesssim \sum_{j=-\frac{e^{-2\pi^2}}{4\pi (1-r)}}^{\frac{e^{-2\pi^2}}{4\pi (1-r)}} \int_{I_j}\int_{I_j} |g(x)g(z)| \sup_{\xi\in \T} |(\psi_{N(x)}-\psi_{N{z}})'(\xi)| dxdz
$$
$$
\lesssim  \sum_{j=-\frac{e^{-2\pi^2}}{4\pi (1-r)}}^{\frac{e^{-2\pi^2}}{4\pi (1-r)}} \frac{1}{1-r} \left( \int_{I_j} |g(x)| dx \right)^{2} \lesssim \|g\|_2^2 \, .
$$
For $x,z\in E$, $k(x)=k(z)$ holds if and only if $k(\tau(x))=k(\tau(z))$ if and only if $k(\eta(x)) = k(\eta (z))$. Hence, $E_0$ is also the diagonal for the other three kernels and similar estimates hold for the other three terms integrated over $E_0$. 

If $(x,z) \in E\times E\setminus E_0$, then by the choice of the set $E$, both $k(x)$ and $k(\eta(x))$ are even, and $|k(x)-k(z)| \geq 2$ and $|k(\eta (x)) - k(\eta (z))|\geq 2$. Thus, using parts d and e of Lemma \ref{taueta}, we can change the Hilbert kernels in \eqref{Sigma kernel} to be all $\frac{1}{x-z}$. Denoting
$$
D(x,z) = \chi_N(x,z) + \chi_{\ti N} (\eta(x),\eta(z)) + \chi_M (\tau(x),\tau(z)) + \chi_{\ti M} (\eta (\tau(x)), \eta(\tau (z)) ) \, ,
$$
we have
$$
\left| \Sigma  - \int_{E\times E\setminus E_0} g(x)g(z) D(x,z) \frac{dxdz}{x-z} \right| \lesssim \|g\|_2^2 + 
\int_{E\times E\setminus E_0} |g(x)g(z)| \frac{1-r}{|x-z|^2}dxdz
$$
$$
\lesssim \|g\|_2^2 + \int_E |g(x)Mg(x)| \lesssim \|g\|_2^2 \, .
$$
Hence, \eqref{claim2 estimate} is reduced to showing
\begin{equation}\label{sum est kern D}
    \Sigma' := \left| \int_{E\times E\setminus E_0} g(x)g(z) D(x,z) \frac{dxdz}{x-z}\right| \lesssim \|g\|_2^2 \, .
\end{equation}
We will proceed by partitioning the integration region $E\times E\setminus E_0$ so that we separate the Hilbert kernel into several scales. Let
$$
E_1 := \{ (x,z)\in E\times E \setminus E_0 \,:\, 2\leq |k(x)-k(z)| < \frac{\min (p(x),p(z))}{4e^{2\pi ^2}} \}\, ,
$$
$$
E_2 := \{ (x,z)\in E\times E\setminus E_0 : \frac{p(x)}{e^{2\pi ^2}}  <  |k(x)-k(z)| \leq 10 p(x) \text{ or }$$
$$\frac{p(z)}{4e^{2\pi ^2}}  \leq  |k(x)-k(z)| < 10 p(z)  \}  ,
$$
$$
E_3 := \{(x,z)\in E\times E\setminus E_0 \, : \, 10 p(x) \leq |k(x)-k(z)| < \frac{p(z)}{4e^{2\pi ^2}} \} \, ,$$
and
$$
E_4 := \{(x,z)\in E\times E\setminus (E_0\cup E_2) \, : \, 10 p(z) \leq |k(x)-k(z)|\} \, .
$$
Let, for $j=1,2,3,4$,
$$
\Sigma_j' := \left| \int_{E_j} g(x)g(z) D(x,z) \frac{dxdz}{x-z}\right| \, .
$$
The next four subsections bound each $\Sigma_j'$.

\subsubsection{Estimate for $\Sigma_1'$}
Let $(x,z)\in E_1$. Let $y\in [-\pi, \pi]$ be between $x$ and $z$. By Lemma \ref{psi phase estimate},
$$
|(\psi_{N(z)} - \psi_{N(x)} )'(y) | \lesssim \frac{1}{(1-r )\big(1 + \dist(\frac{y}{2\pi(1-r)}, [N(x),N(z)]) \big)}
$$
$$\lesssim \frac{1}{(1-r) (1+\min (p(x), p(z)))} \, .
$$
Then, by the mean value theorem,
$$
|\chi_N(x,z) -\chi_N(z,x)| \lesssim\frac{|x-z|}{(1-r)(1+\min (p(x),p(z)) )} \, .
$$
The same estimate holds if we change $N$ to $M, \ti N$ or $\ti M$, and $(x,z)$ to $(\tau(x),\tau(z))$, $(\eta(x),\eta(z))$ or $(\eta(\tau(x)), \eta(\tau(z)))$, correspondingly.

Noting that $E_1$ is symmetric, i.e. $(x,z)\in E_1$ imples $(z,x)\in E_1$, we can rewrite $\Sigma'$ in a symmetric form, and estimate
$$
\Sigma'_1 = \frac{1}{2}\int_{E_1} |g(x)g(z)| \frac{|D(x,z)-D(z,x)|}{|x-z|} dxdz 
$$
$$\lesssim \int_{E_1} |g(x)g(z)| \frac{dzdx}{(1-r)( 1+\min (p(x), p(z)))}
$$
$$
= 2 \int_{E_1, \, p(x)>p(z)} |g(x)g(z)| \frac{dzdx}{(1-r)(1+p(z))} $$
$$
\leq 2 \int_{E} |g(z)| \left( \int_{x: |k(x)-k(z)|<p(z)} |g(x)| \frac{dx}{(1-r)(1+p(z))} \right) dz 
$$
$$
\lesssim \int_\T |g(z)| Mg(z) dx \lesssim \|g\|_2^2 \, .
$$

\subsubsection{Estimate for $\Sigma_2'$}
This is trivial.
$$
\Sigma_2' \lesssim \int_{E_2} |g(x)g(z)| \frac{dxdz}{|x-z|} \lesssim \int_E |g(x)| \frac{1}{p(x)} \int_{z: |z-x|\sim p(x)} |g(z)|dz
$$
$$
\lesssim \int_E |g(x)Mg(x)|dx \lesssim \|g\|_2^2 \, .
$$

\subsubsection{Estimate for $\Sigma_3'$}
Let $(x,z) \in E_3$, then we have
$$N(z)-N(x)=k(z)-k(x)+p(z)-p(x) > 0$$
and 
$$M(\tau(z))-M(\tau(x))=k(z)-k(x)-p(z)+p(x) <0$$
have opposite signs, hence,
$$
|\chi_N(x,z) + \chi_M (x,z)| =| \sgn(N(z)-N(x))\sin \big( (\psi_{N(z)}-\psi_{N(x)})(x) \big) 
$$
$$+ \sgn( M(z)-M(x) )\sin \big( (\psi_{M(\tau(z))}-\psi_{M(\tau(x))})(\tau(x))|
$$
\begin{equation}\label{E2 kernel est}
    \lesssim \sin\left( \frac{1}{2}\left| \big( \psi_{N(z)} - \psi_{N(x)} \big) (x) - \big( \psi_{M(\tau(z))} -\psi_{M(\tau(x) )} \big) (\tau(x) ) \right|\right) \, .
\end{equation}

We know that
\begin{equation}\label{tau oddness property}
    \Psi_{a_j} (\tau(x)) = - \Psi_{a_{2k(x)+1-j}} (x) \, .  
\end{equation}
Using \eqref{tau oddness property} we can estimate the difference of the phases in \eqref{E2 kernel est}, we have
$$
\left| \big( \psi_{N(z)} - \psi_{N(x)} \big) (x) - \big( \psi_{M(\tau(z))} -\psi_{M(\tau(x) )} \big) (\tau(x) ) \right| 
$$
$$= \left|\sum_{j=k(x)+p(x)+1}^{k(z)+p(z)} \Psi_{a_j}(x) + \sum_{j=k(z)-p(z)+1}^{k(x)-p(x)} \Psi_{a_j}(\tau(x)) \right|
$$
$$
=\left|\sum_{j=k(x)+p(x)+1}^{k(z)+p(z)} \Psi_{a_j}(x) - \sum_{j=k(z)-p(z)+1}^{k(x)-p(x)} \Psi_{a_{2k(x)+1-j}}(x) \right|
$$
$$
=\left|\sum_{j=1}^{k(z)+p(z)-k(x)-p(x)} \Psi_{a_{j+k(x)+p(x)}}(x) - \sum_{j=1}^{k(x)-k(z)+p(z)-p(x)} \Psi_{a_{k(x)+p(x)+j}}(x) \right|
$$
$$
=   \left|\sum_{j=p(z)-p(x) - |k(z)-k(x)|+1}^{p(z)-p(x)+|k(z)-k(x)|} \Psi_{a_{j+k(x)+p(x)}}(x) \right| \, ,
$$
then, by Lemma \ref{singlephaseaction}, we continue,
$$
= 2\pi|k(z)-k(x)| - \sum_{j=p(z)-p(x)) - |k(z)-k(x)|}^{p(z)-p(x)+|k(z)-k(x)|} O\left( \frac{1}{p(x)+j} \right) 
$$
$$2\pi|k(z)-k(x)| - O\left( \frac{k(x)-k(z)}{p(z)}\right) \, ,
$$
All in all, we have
$$
|\chi_N(x,z) + \chi_M (x,z)| \lesssim \frac{|x-z|}{(1-r)p(z)} \, .
$$

For $\ti N,\ti M$, we have
$$
10 p_{\ti N}(\eta(x)) = 10 p_N(x)\leq  \frac{|k(\eta(x)) - k(\eta(z))|}{e^{2\pi ^2}} \leq \frac{p_N(z)}{4e^{2\pi ^2}} = \frac{p_{\ti N}(\eta(z))}{4e^{2\pi ^2}} \, .
$$
So the estimate applies for the remaining part of the kernel $D$. Namely,
$$
|\chi_{\ti N}(\eta(x),\eta(z)) + \chi_{\ti M} (\tau(\eta(x)),\tau(\eta(z)) )| \lesssim \frac{|x-z|}{(1-r)p(z)} \, .
$$

Putting the two estimate together we conclude,
$$
\Sigma_3' \lesssim \int_{E_3} |g(x)g(z)| \frac{1}{(1-r)p(z)} dxdz
$$
$$
=\int_{E} |g(z)| \frac{1}{(1-r)p(z)} \int_{x: |x-z| < (1-r)p(z)} |g(x)| dx dz
$$
$$
\lesssim \int_E |g(z)Mg(z)|dz \lesssim \|g\|_2^2 \, .
$$

\subsubsection{Estimate for $\Sigma_4'$}
Let $(x,z) \in E_4$, we will estimate
\begin{equation}\label{sum of N ti N kernels}
    |\chi_N (x,z) + \chi_{\ti N} (\eta (x),\eta (z))| \, .    
\end{equation}
Note, that
$$
N(z)-N(x)= k(z)-k(x) - p(x) +p (z)
$$
and
$$
\ti N(\eta (z)) - \ti N(\eta(x)) = \ti k(z) -\ti k(x) - p(x) + p (z) 
$$
have the same sign.
Thus, we write
\begin{equation}\label{phase diff by cos}
    \eqref{sum of N ti N kernels} \lesssim |\cos \frac{1}{2}\left( (\psi_{N(z)}-\psi_{N(x)} )(x) - (\psi_{\ti N(\eta(z))} -\psi_{\ti N(\eta(x))} )(\eta(x)) \right) | \, .
\end{equation}
We want to estimate the difference of the phases above. By rotation symmetry,
$$
\left( \psi_{\tilde{N}(\eta (z))} - \psi_{\tilde{N} (\eta (x)) } \right) (\eta (x)) = \left( \psi_{\ti k(z)+p(z)-\ti k(x) + k(x)} - \psi_{k(x)+p(x)} \right) (x) \, .
$$
Thus,
\begin{equation}\label{phase diff ti N and N}
    \left( \psi_{\tilde{N}(\eta (z))} - \psi_{\tilde{N} (\eta (x)) } \right) (\eta (x)) - \left( \psi_{N(z)}  -\psi_{N(x)} \right) (x)    
\end{equation}
$$
= \left( \psi_{\ti k(z) +p(z)-\ti k(x)+k(x)} - \psi_{k(z)+p(z)} \right) (x) \, ,
$$
plugging in the asymptotic estimate of Lemma \ref{phaseasymptotics}, we continue
$$
= \sum\limits_{j=k(z)+p(z)+1}^{\ti k(z)+p(z)-\ti k(x) +k(x)} \left(\pi - \frac{1}{2\pi(j-k(x))} + O\left( \frac{1}{|j-k(x)|^2} \right) + O\left(\frac{1-r}{|j-k(x)|} \right) \right)
$$
$$
= (\ti k(z) - \ti k(x) +k(z)-k(x) ) \pi + O(1-r)
$$
$$
+O \left( \frac{1+p(z)}{|k(z)-k(x)|} \right)
- \sum_{j=k(z)-k(x)+1}^{\ti k(z) -\ti k(x) } \frac{1}{2\pi j}
$$
$$
= 2 L \pi + O(1-r)+O\left(\frac{1+p(z)}{|k(x)-k(z)|} \right) - \frac{1}{2\pi} \log \frac{\ti k(z)-\ti k(x)}{k(z)-k(x)}.
$$
Here, $2L = \ti k(z) - \ti k(x)+ k(z)-k(x)$ is an even integer by the definition of $\ti k(x)$ and of the set $E$. Again, by the definition of $\ti k(\cdot)$,
$$
\frac{1}{2\pi} \log \frac{\ti k(z)-\ti k(x)}{k(z)-k(x)} = \frac{1}{2\pi} \log \left(e^{2\pi^2} + O\big(\frac{1}{|k(z)-k(x)|}\big) \right) = \pi + O\left(\frac{1}{|k(z)-k(x)|}\right).
$$
All in all,
\begin{equation*}
    \eqref{phase diff ti N and N} = 2L\pi - \pi + O(1-r) + O\left( \frac{1+p(z)}{|k(z)-k(x)|}\right) \, .
\end{equation*}
Hence, continuing from \eqref{phase diff by cos}, we obtain
$$
|\chi_N (x,z) + \chi_{\ti N} (\eta(x),\eta(z)) | \lesssim |\cos \left( - \frac{\pi}{2} + O(1-r) + O \left( \frac{1+p(z)}{|k(x)-k(z)|} \right) \right) |
$$
$$
\lesssim 1-r + \frac{1+p(z)}{|k(x)-k(z)|} \, .
$$
Similar estimate holds for $M$ and $\ti M$, so that,
$$
\Sigma_4' \lesssim \int_{E_4} |g(x)g(z)| \left(\frac{1-r}{|x-z|} + \frac{(1-r)p(z)}{|x-z|^2} \right) dxdz
$$
$$
\lesssim \int_{E} |g(z)| \sum_{j=0}^{|\log ((1-r)p(z))|} \left( (1-r)2^{j} + (1-r)p(z)2^{2j} \right)\int_{|x-z|\sim 2^{-j}} |g(x)|dx
$$
$$
\lesssim \int_E |g(z)| Mg(z) \sum_{j=0}^{|\log ((1-r)p(z))|} \left( (1-r) + (1-r)p(z)2^{j} \right) \lesssim \|g\|_2^2 \, .
$$

\section{Proof of Corollary \ref{aeconvergence}}\label{proofofcor}
For $n\in \mathbb{N}$, let $T_n$ denote the $n$th MT partial sum operator corresponding to the sequence $\mathbf{b}$. First, we will show that lacunary maximal operator $\sup_n|T_{2^n}|$ is bounded. One approach could be to directly apply \cite[Theorem 3]{MNATSAKANYAN2022109461}. But this would be an overkill, as the latter is technically equivalent to the polynomial Carleson theorem \cite{Lie3, Pav}. Instead, we choose to reduce the matters to the lacunary maximal operator of the classical Fourier series whose boundedness relies only on the Littlewood-Paley theory \cite[Corollary 5.3.3]{Grafakos}.

We have
$$
\Psi_w''(y) = -2|w| \frac{(1-|w|^2)\sin (y-\arg w)}{(1+|w|^2-2|w|\cos (y-\arg w))^2} \, .
$$
Let $r_m=1-2^{-m}$.We compute
$$
\left|\sum_{j=1}^{2^m} \Psi_{b_{2^m+j}}'' (y) \right| \lesssim 2^{-m} \sum_{j=1}^{2^m} \frac{|\sin (y-2\pi j 2^{-m})|}{(1+r_m^2-2r_m\cos (y-2\pi j 2^{-m}) )^2}
$$
$$
\lesssim 2^{-m}\sum_{j=1}^{2^m} \frac{j2^{-m}}{2^{-2m}(1+j^2)} \lesssim m \, .
$$
By the mean value theorem and the rotation symmetry
$$
|\sum_{j=1}^{2^m} \Psi_{b_{2^m+j}}' (y) - \sum_{j=1}^{2^m} \Psi_{b_{2^m+j}}' (0) | \leq 2^{-m} \sup_\T |\sum_{j=1}^{2^m} \Psi_{b_{2^m+j}}'' (\xi)| \leq \frac{m}{2^m} \, .
$$
By yet another mean value theorem
$$
|\psi_{2^m} (x) - \psi_{2^m} (y) - (x-y)\psi_{2^m}'(0)| = |\int_x^y (\psi_{2^m}'(\xi) -\psi_{2^m}'(0)) d\xi |
$$
$$
\leq |x-y| \sum_{k=1}^m \frac{k}{2^k}\lesssim |x-y| \, .
$$
On the other hand, we can compute
$$
\psi_{2^m}'(0) = \sum_{j=1}^{2^m} \Psi_{b_{2^m+j}}' (0) = \sum_{j=1}^{2^m} \frac{1-r_m^2}{1+ r_m^2-2r_m\cos (2\pi j2^{-m})}
$$
$$
=\sum_{j=1}^{2^m} \left( \frac{1-r_m^2}{1+r_m^2-2r_m\cos (2\pi j2^{-m})} - \frac{2^m}{2\pi}\int\limits_{2\pi j2^{-m}}^{2\pi (j+1)2^{-m}} \frac{(1-r_m^2) d\theta}{1+r_m^2-2r_m\cos \theta } \right)
$$
$$
+\frac{2^m}{2\pi}\int_\T \frac{(1-r_m^2) d\theta}{1+r_m^2-2r_m\cos \theta }= 2^m + 
$$
$$
\frac{(1-r_m^2)2^m}{2\pi} \sum_{j=1}^{2^m} \int\limits_{2\pi j2^{-m}}^{2\pi (j+1)2^{-m}} \frac{2r_m (\cos (2\pi j2^{-m}) - \cos \theta) d\theta}{(1+r_m^2-2r_m\cos \theta) (1+r_m^2-2r_m\cos (2\pi j2^{-m}) )}
$$
$$
=2^m+ 2^{-m} O\left(\sum_{j=1}^{2^m} \frac{j2^{-m}}{2^{-2m}+j^2 2^{-2m}} \right)
=2^m+ O(m) \, .
$$

We conclude that the boundedness of the operator $\sup_m |T_{2^m}f|$ is equivalent to that of $\sup_m |S_{\lambda_m}f|$, where $\lambda_m=2^m+O(m)$ is some lacunary sequence and $S_n$ is the $n$th partial sum of the classical Fourier series. For the latter the $L^2$ boundedness is much easier to obtain than for the full Carleson operator and can be found, for example, in \cite[Corollary 5.3.3]{Grafakos}.

To combine the estimate for the lacunary operator with Theorem \ref{main}, let us consider the following chain of inequalities.
$$
\sup_n |T_n f| = \sup_m \sup_{1\leq j\leq 2^m} |T_{2^m+j} f|
$$
$$
\leq \sup_m \sup_{1\leq j\leq 2^m} |T_{2^m+j} f -T_{2^m} f| + \sup_m |T_{2^m}f|
$$
$$
\leq \left( \sum_{m=1}^\infty \sup_{1\leq j\leq 2^m} |(T_{2^m+j} -T_{2^m} ) f| ^2 \right)^{\frac{1}{2}} + \sup_m |T_{2^m}f| \, .
$$
The second summand is bounded in $L^2$ by the above discussion. Let $H_m$ be the subspace of $H^2 (\T)$ spanned by $\phi_{2^m+j}$, $j=1,\dots , 2^m$, and let $P_m$ denote the orthogonal projection from $L^2$ to $H_m$. Then, recalling that $T_n$ was the MT partial sum operator,
$$
(T_{2^m+j}-T_{2^m}) f = (T_{2^m+j}-T_{2^m}) (P_m f) \, .
$$
Furthermore, as the sequence $\mathbf{b}$ satisfies \eqref{completeness condition}, we have
$$H^2(\T) = \bigoplus_{m=0}^{\infty}H_m \, .$$
Applying Theorem \ref{main}, we finish the proof
$$
\|( \sum_{m=1}^\infty \sup_{1\leq j\leq 2^m} |(T_{2^m+j} -T_{2^m} )f| ^2 )^{\frac{1}{2}} \|_2
= \left(\sum_{m=1}^\infty \|\sup_{1\leq j\leq 2^m} |(T_{2^m+j}-T_{2^m} )(P_mf) \|_2^2 \right)^{\frac{1}{2}}
$$
$$
\lesssim \left( \sum_{m=1}^\infty \|P_mf \|_2^2 \right)^{\frac{1}{2}} = \|f\|_2 \, .
$$

\section{A simple proof of \eqref{counterexample}}\label{simple counterexample}
In this section $\phi_j$, $j=0,1,\dots$, will denote the MT system corresponding to the sequence $\mathbf{d_r}$. Our main claim is the following.
\begin{claim}\label{pointwise lower 1.9}
Let $M = \frac{1}{2(1-r)\log \frac{1}{1-r}}$. For $0\leq n \leq M$ and $\frac{2n+2}{2M} \leq y \leq \frac{2n+3}{2M}$, we have
\begin{equation}\label{pointwise lower}
    \Im \left( \sum_{j=0}^{n} \phi_{2j} (e^{iy}) \right) \gtrsim \frac{1}{\sqrt{1-r}} \, .
\end{equation}
\end{claim}
Before proving the claim, let us see that it implies \eqref{counterexample}. Let
$$f(e^{iy}) = \sum_{j=1}^{M} \phi_{2j} (e^{iy}) \, .$$
Then,
$$
T_{2n} (f) = \sum_{j=0}^{n} \phi_{2j} \text{ and } \| f\|_{L^2(\T)}^2 =M \, .
$$
Furthermore, \eqref{pointwise lower} implies
$$
\| \sup_n T_n(f) \|_{L^2(\T)}^2 \geq \sum_{n=1}^M \| T_{2n}(f) \|_{L^2([(2n+2)/2M,(2n+3)/M])}^2
\gtrsim \sum_{n=1}^M \frac{1}{(1-r)M} = \frac{1}{1-r}\, ,
$$
so that
$$\| \sup_n T_n(f) \|_{L^2(\T)}^2 / \|f\|_{L^2(\T)}^2 \gtrsim \log \frac{1}{1-r} \, .$$

\begin{proof}[Proof of Claim \ref{pointwise lower 1.9}]
Fix $n$ and $y$ as in the hypothesis of the claim and let us look at the phase of the Blaschke product $B_{2j}$ for $1\leq j\leq n$. By Lemma \ref{phaseasymptotics},
$$
\psi_{2j} (y) = \sum_{m=1}^{2j} \Psi_{a_m} (y) 
= 2\pi j - \sum_{m=1}^{2j} \frac{1}{1+[y/1-r]-m\log \frac{1}{1-r}} + O \left( \frac{1}{\log \frac{1}{1-r}} \right)
$$
$$
= 2\pi j - \sum_{m=1}^{2j} \frac{1}{\log \frac{1}{1-r} (2n+1-m)} + O \left( \frac{1}{\log \frac{1}{1-r}} \right)
$$
$$
= 2\pi j - \frac{A(n,m)}{\log \frac{1}{1-r}} + O \left(\frac{1}{\log\frac{1}{1-r}} \right) \, ,    
$$
where
$$0\leq A(n,m) = \sum_{m=1}^{2j} \frac{1}{2n+1-m} \leq \log n + O(1) \leq \log \frac{1}{1-r} + O(1),$$
as $n\leq M$. In particular,
\begin{equation}\label{Blaschke phase}
    |\psi_{2j} (y) - 2\pi j| \leq 1< \frac{\pi}{2}.
\end{equation}
On the other hand, for the orthogonalizing factor, we have
\begin{equation}\label{real part of factor}
    \Re \frac{\sqrt{1-r^2}}{1-re^{-i(2j+1)/M}e^{iy}} = \frac{\sqrt{1-r^2}(1-r\cos(y-\frac{2j+1}{M}))}{(1-r\cos (y-\frac{2j+1}{M}))^2 + r^2 \sin^2 (y-\frac{2j+1}{M})} \sim \sqrt{1-r}    
\end{equation}
and
$$
\Im \frac{\sqrt{1-r^2}}{1-re^{-i(2j+1)/M}e^{iy}} = \frac{\sqrt{1-r^2}r\sin (y-\frac{2j+1}{M})}{(1-r\cos (y-\frac{2j+1}{M}))^2 + r^2 \sin^2 (y-\frac{2j+1}{M})}
$$
\begin{equation}\label{im part of factor}
    \gtrsim \sqrt{1-r} \frac{M}{n+1-j} = \frac{1}{(n+1-j)\sqrt{1-r}\log \frac{1}{1-r}} \, .
\end{equation}
Putting together \eqref{Blaschke phase}, \eqref{real part of factor} and \eqref{im part of factor}, we write
$$
\Im \phi_{2j} (y) = \Re e^{i\psi_{2j} (y)} \Im \frac{\sqrt{1-r^2}}{1-re^{-i(2j+1)/M + iy}} + \Im e^{i\psi_{2j} (y)} \Re \frac{\sqrt{1-r^2}}{1-re^{-i(2j+1)/M + iy}}
$$
$$
\gtrsim \frac{1}{(n+1-j)\sqrt{1-r}\log \frac{1}{1-r}} + O(\sqrt{1-r}) \gtrsim \frac{1}{(n+1-m)\sqrt{1-r}\log \frac{1}{1-r}} \, ,
$$
as $0\leq n+1-j \leq M= \frac{1}{(1-r)\log \frac{1}{1-r}}$. Summing up for $1\leq j\leq n$,
$$
\Im \sum_{j=0}^{n} \phi_{2j} (y) \gtrsim \sum_{j=1}^n \frac{1}{(n+1-j)\sqrt{1-r}\log \frac{1}{1-r}} \sim \frac{1}{\sqrt{1-r}} \, .
$$
\end{proof}

\bibliographystyle{amsplain}

\bibliography{references}

\end{document}